\documentclass[runningheads]{llncs}
\usepackage{latexsym,amsfonts,amsmath}
\usepackage{graphicx}
\usepackage{epstopdf}

\begin{document}

\title{A Lie systems approach for the first passage-time of piecewise deterministic processes}
\titlerunning{Lie systems in Mathematical Finance}
\authorrunning{Avram, F., Cari\~nena, J.F. and de Lucas, J.}

\author{Florin Avram\inst{1}, Javier de Lucas\thanks{Corresponding author}\inst{2},  \and  Jos\'e F. Cari\~nena\inst{3}}

\institute{Universit\'e de  Pau, France 64000\\
\email{Florin.Avram@univ-Pau.fr} \\ \
\and Institute of Mathematics, Polish Academy of Sciences, ul. \'Sniadeckich 8, P.O. Box 21, 00-956 Warszawa, Poland\\ \email{delucas@impan.pl}\\ \
\and Departamento de F\'isica Te\'orica, Universidad de Zaragoza, C. Pedro Cerbuna 12, 50009 Zaragoza, Spain\\
\email{jfc@unizar.es}
}

\maketitle

\begin{abstract}
Our paper illustrates how the theory of Lie systems
allows recovering known results and provide new examples of piecewise deterministic processes
 with phase-type jumps for which the
corresponding first-time passage problems may be solved explicitly.
\end{abstract}

\section{Introduction}

{\bf Motivation.} Probabilists share with physicists the tradition of paying special attention to analytically tractable models, like for example solvable Markovian semigroups/evolution equations, in their respective parallel terminologies. Such models have lead to significant progress over the centuries, as witnessed for example by the recent cases of the Black-Scholes and affine models, in mathematical finance.

With a few exceptions, the main approach used nowadays for dealing with Kolmogorov evolution equations (to adopt a `unified' name) in applied probability is the transform method initiated by Laplace, Poincar\'e and Weierstrass.  It is intriguing to ask whether the alternative algebraic-geometric approach pioneered by Lie \cite{FLLS}, which is going these days through a revival period (see, for instance, \cite{FLCGM00,FLCL10Int,FLIb99,FLLazOrt,FLPW}) stimulated by today's symbolic computing, could turn out to be as useful in probability.

Note that Lie's approach would seem to be {\it taylor made} for jump-diffusion
processes, which are generated by combining noncommuting operators.

{\bf Jump-diffusions.}
A recurring theme in applied probability  is distinguishing between two possible sources of uncertainty:
 small continuous changes modeled by diffusion, and `catastrophic' changes modeled by a jump process.
To resolve this issue, one uses jump-diffusions models, i.e. solutions of a SDE (stochastic differential equation)
\begin{eqnarray}  \label{FLjd}
d X_t&= & \varphi(X_t) d t + \sigma(X_t) d B_t - d S_t,
\end{eqnarray}
where
\begin{itemize}
\item  $B_t$ is standard Brownian motion,
\item  $S_t$ is a pure jump process, with a {\it Levy density} $\nu(x, z):=\lambda(x) \,b(z),$ which may arise, for example, from i.i.d. jumps $C_i$ (sometimes one sided, for example negative), whose density,
distribution, complementary distribution, and first moment are denoted respectively by $b(x), B(x),  \bar{B}(x), b_1.$
\end{itemize}

The first two terms (the drift $\varphi$ and the variance $\sigma$) of the Levy-Khinchine
 triple $\varphi, \sigma,\nu$ define a continuous diffusion process, and the last term defines a pure
jump/convolution  process.

Jump-diffusions (\ref{FLjd}) are Markovian processes with associated evolution/ backward Kolmogorov equation
\begin{equation} \label{FLevo}
\frac{\partial f(x,t)}{\partial t} =\mathcal{G}_x f(x,t), \quad
f(x,0)=f_0(x),\end{equation} which describes {\it expectations
 evolving in time}
 $$f(t,x)=\mathbb{E}_{X_0=x}  f_0(X_t).$$
The infinitesimal generator operator is given by
\begin{multline*}
\mathcal{G} f(x)= \mathcal{G}_x f(x)= \varphi(x) f'(x) + \frac{\sigma^2(x)}{2}
f''(x) +\\ +\int_{-\infty}^\infty (f(x - z)-f(x)) \nu(x,z) dz \mathcal{G}^{(d)} f(x)
+ \mathcal{G}^{(j)} f(x),
\end{multline*} for any twice continuously differentiable and bounded function $f(x)$, where the second part $\mathcal{G}^{(j)}_x$ is associated to the pure jump convolution part.

One important example, already studied in Kolmogorov's founding paper \cite{FLKolmogorov}, is that of {\it hypergeometric diffusions} with quadratic variance and linear drift:
\begin{equation*}
\mathcal{G}^{(d)} f(x):= \left(a_{2}x^{2} + a_{1}x + a_{0} \right)
\displaystyle\frac{\partial^{2} f}{\partial x^{2}} + \left(\varphi_{1}x + \varphi_{0}\right) \frac{\partial f}{\partial x}.
\end{equation*}

The Levy model is obtained when the  variance and drift rates as well as the  Levy intensity $\nu(x,z)$ are independent of $x$.

\begin{example} The Cram\'er Lundberg risk model (1903) is one of the most studied models in applied probability \cite{FLLun}, describes the surplus of an insurance company:
\begin{equation*}
U(t)=u + c \; t- S(t):=u + c \; t- \sum_{i=1}^{N(t)}C_i,
\end{equation*}
with initial capita $u$,  linear premium rate/drift $c \, t$ and
{\it claims} $C_k,$ modeled by a sequence of i.i.d. positive random
variables with a common density  $f(x)=f_C(x),$ and which arrive  at
the increase points $N=\{N_t, t\ge0\}$ of an independent counting
process with $\mathbb{E} N_t=\lambda t$.  If moreover $N_t$ is a
Poisson process with exponential interarrival times, then $S(t)$ is
a compound Poisson process with positive summands, and $U(t)$, with
$t \in \mathbb{R}_+$, is Markovian.

If, more generally,  $S(t)$ is a subordinator (nondecreasing Levy process), then $U(t) $ is a {\it spectrally negative} Levy process (with jump measure supported on the negative half-line).
\end{example}

{\bf Phase-type distributions.} For a general line of attack on jump diffusions, it seems natural to start by adding an independent compound Poisson process with phase-type jumps  $C_i$ (the convolution
 term in the evolution equation),  with distribution
$$
\bar{B}(x):=P[C >x]=\mbox{\boldmath $\beta$} e^{\bf B x} \mbox{\boldmath $1$},
$$
where ${\bf B }, \mbox{\boldmath $\beta$},
\mbox{\boldmath $1$}$ are respectively a subgenerating matrix
(nonnegative off-diagonal elements and nonpositive row sums), a
probability row vector and a column vector of ones. Note that the
density (the negative of the derivative of $ \bar{B}(x)$) is
$$b(x)=\mbox{\boldmath $\beta$} e^{\bf B x} {\bf b},$$
where ${\bf b}=- {\bf B} \, \mbox{\boldmath $1$}$ is a column vector
with $n$ components.

Phase-type distributions have the advantage that the
integro-differential equations satisfied for example by the ruin
probabilities $\Psi(t,x)$ may be transformed into ordinary
differential equations --see for example Paulsen \cite{FLPau}, which
resolves by this approach four notoriously difficult particular
cases, in terms of special hypergeometric functions.  Paulsen's
paper arose the natural question of whether under more general
models of hypergeometric diffusions with phase-type jumps, ruin
probabilities could be determined analytically as well.

Some further particular examples were studied for example in the papers \cite{FLAU,FLJacJen}. We are continuing here this line of research,  providing a new solvable model.

Note that it is also possible to consider the more general
matrix-exponential distributions, where one  assumes only the
nonnegativity and integrability (to $1$) of the density $b(x)$, but
the advantages are not clear, since phase-type distributions are
already dense \cite{FLA00b}.

{\bf First passage problems.} Denote by
\begin{equation*}
\tau_L^+= \inf\{t \geq 0;\; X_t  > L\},\qquad \tau=\tau_l =
\inf\{t \geq 0;\; X_t < l\},
\end{equation*}
the first passage times
of a stochastic process above/below given levels $L,l$. The latter,
also called {\it ruin time} in the  insurance literature, is one of the
oldest applications of probability and ODE's, introduced by Thiele,
the founder   of the Danish insurance company Hafnia (1872) --see
 www.stats.ox.ac.uk/~steffen/seminars/centertalk.pdf.

{\bf Ruin probabilities.} {The first objects of interest in first
passage theory are the
 {\it finite-time} and {\it ultimate/infinite horizon} ruin probabilities
 $\Psi(t,x)$
 and the related {\it survival} probabilities $\overline{\Psi}(t,x)$
 \begin{eqnarray*}
&&\Psi(t,x)= P_x [\tau  \leq t], \quad \quad \quad \quad \quad \quad \quad \Psi(x)= P_x
[\tau  < \infty],\\
  &&\overline{\Psi}(t,x)=
  P_x [\tau  > t] = 1 - \Psi(t,x),  \quad \overline{\Psi}(x)= P_x
[\tau  = \infty].
 \end{eqnarray*}

For the Markovian case, a first step/infinitesimal analysis shows that the ultimate ruin probabilities are harmonic functions, satisfying
\begin{equation}\label{FLfk}
\mathcal{G} \Psi(u):=\frac{\sigma^2 }{2} \;  \Psi''(u) + c {\Psi}'(u) -\lambda  {\Psi}(u) +
\lambda \int_0^u {\Psi}(u-z) f(z) d z + \lambda \bar{F}(u) =0,
\end{equation}
with $\Psi(u)=1$ and $u \leq 0$.
\begin{note}
Note that the trivial solution of equation (\ref{FLfk}) with $\Psi(u)=1$ (true when
$\lim_{t \to \infty} X(t)=-\infty$) may be  discarded by adding
the restriction
$$ \lim_{t \to \infty} X(t)=\infty \leftrightarrow \lim_{u \to \infty} \Psi(u)=0.$$
\end{note}

The condition for that to hold is the same as for the process reflected at $0$ to be nonergodic; in the Levy case, this is quite simple:
$$\mathbb{E} \frac{X(t)- X(0)}{t}=
c -\lambda \mathbb{E} C_1 >0,$$
(by the classic law of large numbers), but things get more complicated for the generalizations considered here, with the nonergodicity condition involving the invariant measure $\mathop{\mathrm{\mathfrak{p}}}(x)$ of the process, i.e. the nonnegative solution of the adjoint equation $\mathcal{G}^* \mathop{\mathrm{\mathfrak{p}}}(x)=0$.

\begin{note}
In (\ref{FLfk}) and in any problem involving the process under
consideration must appear the same   operator $\mathcal{G} ,$ which is also
the generator of the associated semigroup of transition operators.
Put informally, $\mathcal{G} $ is the `analysis dictionary entry' associated
to a given Markovian process. For example, the finite-time ruin
probabilities must satisfy the nonautonomous {\it backward Kolmogorov}
partial integro-differential equation:
\begin{equation}\label{FLtd}
-\frac{\partial}{\partial t} \Psi(u,t) + \mathcal{G} \Psi(u,t) =0,
\quad \Psi(u,t)=1, u \leq 0, \quad \forall t.
\end{equation}
\end{note}

{\bf Laplace transform in time.} The finite-time ruin probabilities
are very seldom available analytically; one natural way to approximate them is by computing and inverting their Laplace-Carson transform
\begin{equation}
\label{FLLCar}\Psi_q(x)=\int_0^\infty q \mathrm{e}^{-q t} \Psi(t,x) dt,
\end{equation}
which may be also viewed as the probability of passage of the process `killed'
after an independent exponential random variable ${\bf e}_q$ of rate $q>0$:
\begin{multline*}
\Psi_q (x) =
\Pr_x ( \tau < {\bf e}_q) = \mathbb{E}_x (e^{-q \tau} I_{\{\tau < \infty\}})=\\
\int_0^{\infty} e^{-qt} \Pr_x(\tau \in d t)=\int_0^{\infty} q e^{-qt}
\Psi (t,x) dt,
\end{multline*}
where the last equality follows by integration by parts.

More generally, we want to calculate the {\it killed ruin probability}
\begin{equation}
\Psi_q (x) = \Psi_q^l (x)= \mathbb{E}_x e^{-\int_0^ \tau q(X(s)) ds},
\end{equation}
for general discount functions $q(X(s))$, and the killed ruin with overshoot penalty
$$
\Psi_{q}^{l,\xi}(x)=\mathbb{E}_x (e^{-\int_0^\tau q (X(s)) ds + \xi
(X_{\tau}-l)}I_{\{\tau < \infty\}}).
$$

By Dynkin's formula, this function is a solution of the {\it Sturm-Liouville type equation}
\begin{eqnarray} \label{FLSL}
\mathcal{G} f(x) -q(x) f(x)=0,\\
\quad f(x)=e^{\xi (x-l)}, x < l, \quad f(l)=1, \quad &\text{if }
\sigma \neq 0,
\end{eqnarray} see for example \cite[pag. 824]{Cai05}.
Note that the number of boundary conditions depends on the presence
of Brownian motion (probabilistically, this is required by the
presence of a new unknown: the probability of crossing continuously).

\begin{note}
The arguments $l, \xi$ will be often suppressed (in particular when they equal $0$).
\end{note}

The evolution equation (\ref{FLevo}) and the corresponding time-independent counterparts,
the invariant measure and the harmonic functions of interest in first passage theory, have been intensively studied for {\it diffusions} and for {\it Levy processes}.

In both the diffusion  and the Levy case there exists a beautiful first-passage theory. In the latter case, this reduces the calculation of many of the functions (\ref{FLSL}) to the Wiener-Hopf factorization (of the generator, or rather of its symbol) --see for example \cite{FLAPPdiv,FLKbook}.

A question which begs of itself is  combining these two well understood generators, i.e. investigating the existence of extensions under some {\it common umbrella case} like processes with rational generators.  Note that the affine case has already been thoroughly investigated, for example in mathematical finance, for  modeling  interest rates.

Below, we focus on the case of piecewise deterministic processes (no diffusion) with {\it phase-type downward jumps}. In this case, the Feynman-Kac integro-differential equation \cite{FLKac49} for killed
ruin probabilities  may be brought to the form of a ODE linear system \cite{FLA00b} of the form
\begin{equation}\label{FLLine}
\left(
\begin{aligned}
 \Psi'(x)\\
{\bf M}'(x)\\
\end{aligned}\right)=
\left(
\begin{matrix}
\frac{\lambda+q}{\varphi(x)} & \frac{-\lambda {\mbox{\boldmath $\beta$}}}{\varphi(x)}\\
 {\bf b} &  {\bf B}\\
\end{matrix}\right)
\left(
\begin{aligned}
 \Psi(x)\\
{\bf M}(x)\\
\end{aligned}\right),
\end{equation}
where  ${\bf B}$ is a $n\times n$ stochastic generating matrix,
where ${\bf M}$ is a column vectors with $n$ components, and
${\mbox{\boldmath $\beta$}}=(\beta_1,\ldots,\beta_n)$ is a row
probability vector, i.e. $\sum_{i=1}^n\beta_i=1$. The variable
$\Psi$ is the killed ruin probability, the function $\varphi(x)$ is
the {\it drift}, the constant $q$ is the killing rate/Laplace
transform argument, and the components $M_1,\ldots,M_n$ of the
vector ${\bf M}$ are killed ruin probabilities, obtained  by
changing the jumps to segments of slope $\pm 1$ for
upwards/downwards jumps,  and by associating `auxiliary stages of
artificial time' to the phases of the jumps ($M_i(x)$ is thus the
killed ruin probability when starting at $x$ in phase $i$).

Note that with {\it upward jumps}, the analog equation will have  the
last rows corresponding to the $M$ variables multiplied by $-1$.

\begin{note}
We will restrict to  first-passage problems in  domains where the  drift $\varphi(x)$ does not change sign,  which determines then corresponding boundary conditions.
\end{note}

The theory of Lie systems \cite{FLCGM00,FLCL10Int,FLLS,FLPW}, of which nonautonomous systems of first-order linear homogeneous differential equations are a particular case, states that linear systems like
(\ref{FLLine}) are {\it integrable by quadratures}, i.e. their solutions can be explicitly obtained in terms of `algebraic operations' and integrations of certain given functions depending on one variable \cite{FLKozlov}, if there exists a
finitely generated solvable matrix Lie algebra $\mathfrak{g}$ such
that
\begin{equation}\label{FLIntCon}
A_x\equiv \left(
\begin{matrix}
\frac{\lambda+q}{\varphi(x)} & \frac{-\lambda{\mbox{\boldmath $\beta$}}}{\varphi(x)}\\
{\bf b} & {\bf B}\\
\end{matrix}\right)=\frac{\lambda}{\varphi(x)}
\left(
\begin{array}{cc}
\frac{\lambda+q}{\lambda} & -{\mbox{\boldmath $\beta$}}\\
0 & 0\\
\end{array}\right)+\left(\begin{matrix}
0 & \,\,0\\
{\bf b} & \,\,{\bf B}\\
\end{matrix} \right)\in \mathfrak{g}, \qquad \forall x\in \mathbb{R}.
\end{equation}

For example, in the easy case of  a constant  drift $\varphi(x)=c$, the integrability condition (\ref{FLIntCon}) is trivially satisfied as $A_x\in\mathfrak{g}\equiv\langle A_{0}\rangle,$ for all $x\in
\mathbb{R}$, and the solutions are phase-type functions.

In the case of a non-constant drift $\varphi(x)$, verifying condition (\ref{FLIntCon}) reduces to proving that the sequence of commutators
\begin{equation}
\label{FLcoms} [\bar T_1,\bar T_2], [\bar T_1,[\bar T_1,\bar T_2],[\bar
T_2,[\bar T_1,\bar T_2]],[\bar T_1,[\ldots [\bar T_1,\bar
T_2]\ldots]],\ldots
\end{equation}
generated by the matrices
$$
\bar T_1=\left(
\begin{matrix}
\frac{\lambda+q}{\lambda} & -{\mbox{\boldmath $\beta$}}\\
\\
0 & 0\\
\end{matrix}\right),\qquad \bar T_2=\left(
\begin{matrix}
0 &\,\, 0\\
{\bf b} &\,\, {\bf B}\\
\end{matrix}\right),
$$
must contain a finite subset (generators) spanning a finite-dimensional Lie algebra with respect to the matrix commutator.

The example we focus on below is that  of downward exponential jumps of rate $\mu,$ over an exponential horizon
${\bf e}_q,$ when the  linear system (\ref{FLLine}) becomes:
\begin{equation} \label{FLexpoLine}
\left(\begin{array}{c}
   \ {\Psi}'(x)\\   M'(x)
\end{array} \right) = \left(\begin{array}{cc}
\frac{\lambda +q}{\varphi(x)} &- \frac{\lambda }{\varphi(x)} \\  \mu & -\mu
\end{array} \right)  \left(\begin{array}{c}
  \Psi(x)\\   M(x)
\end{array} \right)
\end{equation}
with $\mu>0 $,  $\lambda>0$, $q>0$
(here  ${\mbox{\boldmath $\beta$}}=1$,  the unique `probability vector' of dimension 1).

\begin{example}\label{FLex:Seg}
When $q=0,$ the solution to system (\ref{FLexpoLine}) satisfying the conditions $\Psi(\infty)=M(\infty)=0$ and $M(0)=1$ can be derived by subtracting its equations, yielding:
\begin{equation*}
\Psi(x)-M(x)=(\Psi(0)-1) e^{Z(x)}, \; Z(x)=-\mu
x + \int_0^x \frac{\lambda} {\varphi(v)} dv,
\end{equation*}
and, provided $Z(\infty)=
-\infty$, then
$$
M(x)=\mu  (1-\Psi(0))
\int_x^{\infty} e^{Z(v)} dv\Rightarrow \Psi(x)=(1-\Psi(0))\left(\mu
\int_x^{\infty} e^{Z(v)} dv - e^{Z(x)} \right).
$$

With $\varphi(x)=c$ constant, we have $M(x)= e^{-(1-\eta)
\mu x},$ where $\eta$ is the smallest positive root of $c\eta ( -\mu
+ \mu \eta) -((\lambda +q)\eta  - \lambda) =0$, and  when $q=0$, we
have $\eta_1=\frac{\lambda}{ c \mu }, \eta_2=1$, from where we recover
the well known $\Psi(x)= \frac{\lambda}{c \mu } e^{ (\frac{\lambda}{p}-\mu )
x}.$
\end{example}

We will see below that when $q=0$, the theory of Lie system shows that the above system satisfies the integrability condition (\ref{FLIntCon}), and it can be therefore integrated by quadratures. Otherwise, Theorem \ref{FLqn0} ensures that this is no longer the case.

It is natural to ask whether other families of processes with solvable first passage probabilities exist.  By the alternative Riccati approach carried out in Section \ref{FLs:Ric}, we take advantage of the fact that several families of Riccati equations integrable by quadratures have been accumulated in the literature, like the {\it generalised Allen-Stein family} described in Theorem \ref{FLTU}, Hovy's equation obtained with $\varphi(x)=x,
q=0$ (in fact, our equation is a `nonlinear Hovy's equation'), or Robin's family \cite{FLRobin}. Each of these cases will lead us immediately to a class of piecewise deterministic
processes with exponential jumps which can be solved analytically. We exemplify this for  the generalised Allen-Stein family, formed by those particular cases of the system (\ref{FLexpoLine}) satisfying the integrability condition
\begin{equation*}
{\varphi'}/2+ (\lambda+q)- \mu {\varphi}=\kappa
c_1\sqrt{{-\mu\,\lambda\, c_0c_2}{\varphi}},
\end{equation*}
for certain reals constants $c_0, c_1, c_2$, and with $\kappa={\rm
sg}(\lambda\mu/(\varphi(x)c_0c_2))$.

Note that the unifying general approach based on the theory of Lie
systems recently provided in \cite{FLCL10Int}  should allow
recognizing and solving higher dimensional first
passage problems.

{\bf Contents:} Our paper illustrates how the theory of Lie systems
allows recovering known results and provide new ones for first
passage problems of jump-diffusions with phase-type jumps.

We consider various particular instances of piecewise deterministic processes with
exponential jumps, both  associated to solvable and, more interesting,
to non-solvable Vessiot--Guldberg Lie algebras in Sections \ref{FLs:nonsolv} and \ref{FLs:nonsolv2}.

\section{Fundamentals of Lie systems}

Let us briefly recall the most fundamental notions of the theory of Lie systems and the geometrical treatment of differential equations to be used throughout our work.

A fundamental concept in the geometric study of differential equations and in what follows is the hereby called {\it $x$-dependent vector field}. This concept refers to a map
$$\begin{array}{rccc}
X:&\mathbb{R}\times N&\longrightarrow &TN\\
  &(x,y)&\mapsto& X(x,y)\in T_xN,
\end{array}
$$
such that $X_x:y\in N\mapsto X_x(y)=X(x,y)\in TN$ is a standard vector field on $N$, for every $x\in \mathbb{R}$. Note that, in consequence, giving a $x$-dependent vector field is equivalent to providing a family of vector fields $\{X_x\}_{x\in\mathbb{R}}$ on $N$ parametrized by $x\in\mathbb{R}$, what explains its name.

In similarity to standard vector fields, each $x$-dependent vector field, e.g. $X(x,y)=\sum_{i=1}^nX^i(x,y)\partial/\partial x^i$, also admits integral curves (see \cite{FLCa96}) determined by the solutions of the system
\begin{equation}\label{FLGenLie}
\frac{dy^i}{dx}=X^i(x,y),\qquad i=1,\ldots,n.
\end{equation}
Indeed, the relevance of studying the properties of $x$-dependent vector fields mostly relies on investigating systems of the above form by means of the analysis of the properties of its corresponding $X(x,y)$, its denominated {\it associated vector field}.

The theory of Lie systems mainly deals with a class of systems of first-order differential equations, the so-called {\it Lie systems}, associated with $x$-dependent vector fields of the form
\begin{equation}\label{FLLieSys}
X(x,y)=\sum_{\alpha=1}^rb_\alpha(x)X_\alpha(y),
\end{equation}
where $X_1,\ldots, X_r$ are a set of vector fields on $N$ spanning a finite-dimensional Lie algebra of vector fields $V$, the associated {\it Vessiot-Guldberg Lie algebra}. These vector fields can be, at least locally, considered as the fundamental vector fields of a certain action $\Phi:G\times N\rightarrow N$, where $G$ is a Lie group with Lie algebra isomorphic to $V$. In those cases where the form of $\Phi$ can be explicitly determined, the general solution, $y(x)$, of system (\ref{FLGenLie}) can be cast into the form $y(x)=\Phi(g(x),y_0)$, with $y_0$ being any point of $N$ and $g(x)$ being the solution of the  equation on $G$ given by
\begin{equation}\label{FLeqG}
\frac{dg}{dx}=-\sum_{\alpha=1}^rb_\alpha(x)X_\alpha^R(g),\qquad g(0)=e,
\end{equation}
where the $X_\alpha^R$ are certain right-invariant vector fields on $G$, for details see \cite{FLCGM00,FLLazOrt}.

A particular instance of Lie system having a special relevance to our work is given by the following system of homogeneous linear differential equations
\begin{equation}\label{FLLinear}
\frac{dy^i}{dx}=\sum_{j=1}^nA^i_j(x)y^j,\qquad i=1,\ldots,n,
\end{equation}
associated with the $x$-dependent vector field $X(x,y)=\sum_{i,j=1}^nA^i_j(x)y^j\partial/\partial y^i$. In fact, consider the family of vector fields of the form $\{X_{\rm ij}=y^j \partial/\partial y^i\,|\, i, j=1,\ldots,n\}$ closing on the commutation relations $[X_{\rm ij},X_{\rm kl}]=\delta^l_iX_{\rm kj}-\delta^j_kX_{\rm il}$,  for $ i,j,k,l=1\ldots,n$, see \cite{FLCGR99}. The vector fields $X_{\rm ij}$ span a Lie algebra of vector fields $V$, and taking into account that $X_x=\sum_{i,j=1}^nA^i_j(x)X_{\rm ij}$, it follows that (\ref{FLLinear}) is a Lie system related to the Vessiot--Guldberg Lie algebra $V$.

It is important to note that some linear homogeneous systems can be considered, additionally, as a Lie system related to another Vessiot--Guldberg Lie algebra $V_0\subset\subset V$. In order to illustrate this claim, let us consider the linear homogeneous system
\begin{equation}\label{FLSolv}
\left\{
\begin{aligned}
\frac{dy^1}{dx}=&A_1^1(x)y^1+A_2^1(x)y^2,\\
\frac{dy^2}{dx}=&A_1^2(x)y^1-A_1^1(x)y^2,\\
\end{aligned}\right.
\end{equation}
describing the integral curves of the $x$-dependent vector field
$$X_x=A_1^1(x)\left(y^1\frac{\partial}{\partial y^1}-y^2\frac{\partial}{\partial y^2}\right)+A_2^1(x)y^2\frac{\partial}{\partial y^1}+A_1^2(x)y^1\frac{\partial}{\partial y^2}.$$
In this case, taking as $V_0$ the Lie algebra
$$V_0=\left\langle y^1\frac{\partial}{\partial y^1}-y^2\frac{\partial}{\partial y^2},y^2\frac{\partial}{\partial y^1},y^1\frac{\partial}{\partial y^2}\right\rangle,
$$
it turns out that $X_x\in V_0\subset\subset V$ for all $x\in\mathbb{R}$. Consequently, system (\ref{FLSolv}) is also a Lie system related to the Vessiot--Guldberg Lie algebra $V_0$.

There is an alternative way to determine the Vessiot-Guldberg Lie algebras associated with systems of the form (\ref{FLLinear}). In order to do so, let us cast these systems into the matrix form
$$
\left(\begin{array}{c}
\frac{dy^1}{dx}\\
\ldots\\
\frac{dy^n}{dx}
\end{array}\right)=
\left(\begin{array}{ccc}
A^1_1(x)&\ldots&A^1_n(x)\\
\ldots&\ldots&\ldots\\
A^n_1(x)&\ldots&A^n_n(x)
\end{array}\right)
\left(\begin{array}{c}
y^1\\
\ldots\\
y^n
\end{array}\right)={\bf A}(x)\left(\begin{array}{c}
y^1\\
\ldots\\
y^n
\end{array}\right),
$$
and consider the map $\rho:\mathfrak{gl}(n,\mathbb{R})\rightarrow V$ satisfying  that $\rho({\bf M_{ij}})=-X_{\rm ij}$, with $({\bf M_{ij}})^k_l=\delta^k_i\delta^j_l$ and $i,j=1,\ldots,n$. As $[{\bf M_{ij}},{\bf M_{kl}}]=\delta^j_k{\bf M_{\rm il}}-\delta_i^l{\bf M_{\rm kj}}$, for $i,j,k,l=1\ldots,n,$ it follows that the map  $\rho$ is a Lie algebra isomorphism. In consequence, if $X_x\in V_0$ for every $x\in\mathbb{R}$, then ${\bf A}(x)$ belongs to the matrix Lie algebra $\rho^{-1}(V_0)$. Similarly, if ${\bf A}(x)\in \mathfrak{g}$ for every $x\in\mathbb{R}$, then $X_x$ is contained in the Lie algebra $\rho(\mathfrak{g})$. Trivially, the next proposition follows.

\begin{proposition}\label{FLTranslation} System (\ref{FLLinear}) is related to a (solvable) Vessiot--Guldberg Lie algebra $V_0$ if and only if there exists a (solvable) matrix Lie algebra $\mathfrak{g}\subset\mathfrak{gl}(n,\mathbb{R})$ such that ${\bf A}(x)\in\mathfrak{g}$, for every $x\in\mathbb{R}$.
\end{proposition}

The theory of Lie systems establishes that those Lie systems related to solvable Vessiot--Guldberg Lie algebras, and whose associated action $\Phi$ is expressed in terms of elementary functions, are integrable by quadratures. Otherwise, our methods, and in general all the other ones found in the literature, requires the Lie system to hold some kind of extra condition to be still integrable by quadratures. This fact will be used in next sections in order to integrate certain difussion processes (\ref{FLexpoLine}) which cannot by related to solvable Vessiot-Guldberg Lie algebras.

\section{Lie systems and integrability of piecewise deterministic processes\label{FLs:nonsolv}}

Our goal in this section is to study the integrability properties of piecewise deterministic models of the form (\ref{FLexpoLine}). In particular, our aim is to show that the key cases that have been exactly solved in the literature for any drift $\varphi(x)$ (those with $q = 0$ studied by Segerdahl and Paulsen \cite{FLPau,FLSeger}), satisfy an integrability condition described by the theory of Lie systems ensuring that, actually, their solutions can be explicitly written down. Finally, our second objective is to show that if $q\neq0$, the previous integrability condition does not hold and each case determined by a drift $\varphi(x)$ must be analysed separately.

The theory of Lie systems states that integrating a Lie system, e.g. one of the form (\ref{FLGenLie}), can be, at least locally, reduced to solving an equation of the form (\ref{FLeqG}) on a Lie group $G$, provided the action related to the Lie system is known. Moreover, it is also know \cite{FLCGR99} that if some Vessiot-Guldberg Lie algebra related to the Lie system is solvable, the equation (\ref{FLeqG}) can be exactly solved by means of an algorithmic method and, through the action $\Phi$, the general solution to the initial problem can be recovered \cite{FLCGR99}.

Taking into account the above considerations, it is interesting to study when the Lie system (\ref{FLexpoLine}) is associated with a solvable Vessiot--Guldberg Lie algebra and their solutions can be, consequently, obtained. In view of Proposition \ref{FLTranslation}, this fact depends on the existence of a solvable matrix Lie algebra $\mathfrak{g}$ including the  matrices
$$
A_x\equiv \left(
\begin{matrix}
\frac{\lambda+q}{\varphi(x)} & -\frac{\lambda}{\varphi(x)}\\
\mu & -\mu\\
\end{matrix}\right)=\frac{\lambda}{\varphi(x)}\left(
\begin{matrix}
\frac{\lambda+q}{\lambda} & -1\\
0 & 0\\
\end{matrix}\right)+\mu \left(
\begin{matrix}
0 & 0\\
1 & -1\\
\end{matrix}\right),\qquad x\in\mathbb{R}.
$$
Note that when $q=0$ is assumed, the above family reads
$$
A_x\equiv \left(
\begin{matrix}
\frac{\lambda}{\varphi(x)} & -\frac{\lambda}{\varphi(x)}\\
\mu & -\mu\\
\end{matrix}\right)=\frac{\lambda}{\varphi(x)}\left(
\begin{matrix}
1 & -1\\
0 & 0\\
\end{matrix}\right)+ \mu \left(
\begin{matrix}
0 & 0\\
1 & -1\\
\end{matrix}\right)=\frac{\lambda}{\varphi(x)} T_1 + \mu T_2,\qquad
x\in\mathbb{R},
$$
where $$ T_1=\left(
\begin{matrix}
1 & -1\\
0 & 0\\
\end{matrix}\right),\qquad T_2=\left(
\begin{matrix}
0 & 0\\
1 & -1\\
\end{matrix}\right).
$$
The above matrices satisfy $[T_1,T_2]=-T_1-T_2$ and span a two
dimensional solvable Lie algebra $V=\langle T_1,T_2\rangle$.
In consequence, system (\ref{FLexpoLine}), with $q=0$, can be easily solved for every drift $\varphi(x)$. In this way, it is not surprising that the solutions for these models are known since a long time ago \cite{FLPau,FLSeger}.

\begin{proposition} \label{FLqn0}
 When $q \neq 0, $ and for a non-constant drift $\varphi(x)$,
  the matrices $ A_x$  span  the {\bf non-solvable}
Lie algebra $\mathfrak{gl}(2,\mathbb{R})$ of $2 \times 2$ real matrices.
\end{proposition}
\begin{proof}
It is obvious that $A_x=\lambda/\varphi(x)U_1+\mu U_2$, where
$$
U_1=\left(
\begin{array}{cc}
\frac{\lambda+q}{\lambda} & -1 \\
\\
0 & 0\\
\end{array}\right),\qquad U_2=\left(
\begin{array}{cc}
0 &\,\, 0\\
1 &\,\, -1\\
\end{array}\right).
$$
In terms of the above matrices, we define
$$
U_3\equiv([U_1,U_2]+U_2+U_1)\lambda /q+U_2=\left(
\begin{array}{cc}
1 &\,\,\, 0 \\
\\
0 &\,\,\, -1\\
\end{array}\right),\quad
U_4\equiv[U_1,U_3]=\left(
\begin{array}{cc}
0 &\,\,\,\, 2 \\
\\
0 &\,\,\,\, 0\\
\end{array}\right).
$$
Consequently, for every Lie algebra $\mathfrak{g}$ such that $\{A_x\}_{x\in\mathbb{R}}\subset\mathfrak{g}$, the matrices $U_1,U_2,U_3,U_4$ must be contained in $\mathfrak{g}$, as they are made up from Lie brackets and linear combinations of elements of $\mathfrak{g}$. Moreover, as $q\neq 0$, the matrices $U_1,U_2,U_3$ and $U_4$ are linearly independent and they span $\mathfrak{gl}(2,\mathbb{R})$. It follows that $\mathfrak{gl}(2,\mathbb{R})\subset \mathfrak{g}$. Consequently, the Lie algebra $\mathfrak{g}$ is not solvable.
\end{proof}

The consequence of the above proposition is clear: as the Lie algebra $\mathfrak{g}$ is not solvable, there exists no general method to solve system (\ref{FLexpoLine}) for an arbitrary drift when $q\neq 0$. Nevertheless, despite the absence of a general method for solving such systems, the theory of Lie systems also provides criteria to ensure the integration of certain of these systems, generally satisfying some kind of extra condition. The applications of one of these criteria will be the main purpose of the following sections.

It is important to remark that, although there are other alternative methods to obtain the results described within this work, the theory of Lie systems provides a unifying approach to derive and analyse them. Moreover, we think that this theory gives us an appropriate approach to the analysis of systems (\ref{FLLine}), as it can provide generalizations of the methods here described to deal with higher-dimensional cases.

\section{The Riccati approach \label{FLs:Ric}}

In order to simplify our treatment of systems (\ref{FLexpoLine}), we here accomplish an alternative approach to study (\ref{FLexpoLine}) consisting on writing the system
in the coordinate system $\{\eta=\Psi/M,M\}$, bringing it to the
form
\begin{equation}\label{FLRicc}
\left\{
\begin{aligned}
\frac{{\rm d}\eta}{{\rm d}x}&=-\mu\eta^2+\left(\mu+\frac{\lambda+q}{\varphi(x)}\right)\eta-\frac{\lambda}{\varphi(x)},\\
 \frac{{\rm d}M}{{\rm d}x}&=(\eta-1)\mu\, M.
\end{aligned}
\right.
\end{equation}
The above non-linear system is made up from a homogeneous equation
in the variable $M$ and a Riccati equation in the variable $\eta$
(with no dependence on the variable $M$), which  will be called
below {\it Segerdahl's equation}.

After the substitution $y(x)=\mu (\eta(x)-1)$ and the homogenizing
substitution $y(x)=\frac{g'(x)}{g(x)},$ the Riccati equation and it
homogeneous counterpart are brought to the {\it  canonical forms}
\begin{equation}
\label{FLRicstan} y'(x) =- y^2(x) + y(x)\left(\frac{\lambda +q }{\varphi(x)}
-\mu\right) + \frac{q \mu }{\varphi(x)} \Leftrightarrow g''(x) -  z(x) g'(x) -
u(x)g(x)=0,
\end{equation}
where
\begin{equation}\label{FLzdef}
z(x)=\frac{\lambda+q}{\varphi(x)}-\mu,\qquad u(x)=\frac{q \mu
(z(x)+\mu)}{\lambda + q} .
\end{equation}

Note that when $q=0,$ equation (\ref{FLRicstan}) becomes essentially of first order $g''(x) -g'(x)z(x)=0$, and $g'(x)=e^{Z(x)}$, with $\;Z(x)=\int z(x) dx$, recovering Segerdahl's result, see Example \ref{FLex:Seg}.

Note that having an explicit general solution $\eta(x)$ to  the  Riccati equation of the system (\ref{FLRicc}) leads to an explicit general solution  of the  system, obtained by
$$M(x)=L\exp\left(\int^x\eta(t)dt-\mu x\right),
$$
where $L$ is an arbitrary constant. Thus, the solution of the first-passage problem will be available analytically (up to quadratures), whenever the Riccati solution is. In other words, we have reduced the problem of solving system (\ref{FLexpoLine}) to solve a Riccati equation. Additionally, it is worth to remark that, in a similar way, solving higher order dimensional systems (\ref{FLLine}) can be also reduced to solving matrix Riccati equations.

Riccati equations are Lie systems generally related to a Vessiot--Guldberg Lie algebra isomorphic to $\mathfrak{sl}(2,\mathbb{R})$, see \cite{FLCGR99}. It follows that this Lie algebra is not solvable, and the method exposed in \cite{FLCGR99} does not apply. Nevertheless, the theory of Lie systems provides the following result allowing us to integrate Riccati equations for certain specific cases of their $x$-dependent coefficients. This result motivates the use of the alternative approach to systems (\ref{FLexpoLine}) carried out in this Section. Although a similar result could also be demonstrated for these systems, it much simpler to use an alternative approach and use the below result, which proof can be found in \cite{FLCL10Int}.

\begin{theorem}\label{FLTU} The necessary and sufficient condition
for the existence of a  transformation
\begin{equation*}
\bar\eta=G(x)\eta,\quad G(x)>0,
\end{equation*}
 relating the Riccati equation
\begin{equation}\label{FLRiccMod}
\frac{d\eta}{dx}=b_0(x)+b_1(x)\eta+b_2(x)\eta^2\,,  \qquad b_0b_2\ne
0,
\end{equation}
to an integrable one given by
\begin{equation}
\frac{d\bar\eta}{dx}=D(x)(c_0+c_1\bar\eta+c_2\bar\eta^2)\,,\quad c_0c_2\neq
0,\label{FLeqDcs}
\end{equation}
where $c_0, c_1, c_2$ are real numbers and $D(x)$ is a non-vanishing
function, are
\begin{equation}
D^2c_0c_2=b_0b_2,\qquad \left(b_1+\frac{1}{2}\left(\frac{
b_2'}{b_2}-\frac{b_0'}{b_0}\right)\right)\sqrt{\frac{c_0c_2}{b_0b_2}}=\kappa
c_1,\label{FLDinTh2}
\end{equation}
where $\kappa={\rm sg}(D)=sg(b_0/c_0)$. The  transformation is then
uniquely defined by
\begin{equation*}
\bar\eta=\sqrt{\frac{b_2(x)c_0}{b_0(x)c_2}}\,\eta\,.
\end{equation*}
\end{theorem}

As Riccati equations of the form (\ref{FLeqDcs}) can be transformed
easily into an autonomous equation by means of a
$x$-reparametrization and these equations are always integrable by
one quadrature, it follows that all the Riccati equations satisfying
the conditions (\ref{FLDinTh2}) can be solved. The family of Riccati
equations satisfying the previous theorem slightly generalises some
previous results on the topic. Indeed, the similarity with the
result given by Allen and Stein in \cite{FLAS64} is what motivated
calling this family the {\it generalised Allen-Stein family}.

Our main aim now is to determine which Segerdahl's equations can be integrated
by means of the above proposition in order to provide new models admitting
exact solutions. In order to do so, note that Segerdahl's equations can be cast
into the
form (\ref{FLRiccMod}), with
\begin{equation}\label{FLcoeff}
b_0(x)=-\frac{\lambda}{\varphi(x)},\qquad
b_1(x)=\left(\mu+\frac{\lambda+q}{\varphi(x)}\right),\qquad b_2(x)=-\mu.
\end{equation}
Substituting the above functions in the integrability
condition (\ref{FLDinTh2}), we get that Riccati equation (\ref{FLRicc})
is integrable if there exist constants $c_0, c_1,$ and $c_2$ such that the
drift $\varphi(x)$ satisfies the equation
\begin{equation}\label{FLinte}
{\varphi'}/2+ (\lambda+q)+ \mu {\varphi}=\kappa
c_1\sqrt{{-\mu\,\lambda\, c_0c_2}{\varphi}}.
\end{equation}
For example, in the particular case $c_1=0$, the above integrability
condition reads
\begin{equation}\label{FLIntCond}
{\varphi'}+2\mu {\varphi} +2(\lambda+q)=0,
\end{equation}
whose general solution, $\varphi_K(x)$, is
\begin{equation}\label{FLGenSolution}
\varphi_K(x)=\frac{\lambda+q}{\mu}\left(K e^{-2\mu x}-1 \right),
\end{equation}
with $K$ a nonzero real constant.
An explicit solution for the classical ruin problem with this drift follows.

\section{The Ruin probability}\label{FLs:nonsolv2}

We consider now the  ruin problem with downwards exponential jumps and drift (\ref{FLGenSolution}).
In the absence of jumps, we would have a dynamical system an  attractive point $x_0=\frac{ln(K)}{2 \mu}$.

The most interesting exit problem is finding

\begin{equation} \label{FLabove}\Psi(x)=P_x[\tau_L <\tau_l, \tau_L<{\mbox{\boldmath $e$}}_q],\end{equation}
where ${\mbox{\boldmath $e$}}_q$ is an independent exponential
horizon, i.e. $P[{\mbox{\boldmath $e$}}_q >x]=e^{- q x}$, and where
$l<L< x_0.$ Note that in this case the jumps and drift go in
opposite directions and both boundaries may be attained; therefore,
this  requires solving the system (\ref{FLexpoLine}) with the
boundary conditions  $M(l)=0, \Psi(L)=1$.

We turn now to the easiest case $x_0<l<L $, with the drift negative,
in  which case the  problem (\ref{FLabove}) (with $M(l)=0,
\Psi(L)=1$) has no solution,  since because of  the negative drift,
with `no vehicle going up', exit above is impossible.

The other two sided problem \begin{equation} \label{FLbelow}\Psi(x)=
P_x[\tau_l<\tau_L , \tau_l<{\mbox{\boldmath $e$}}_q],\end{equation}
with boundary conditions $M(0)=\Psi(0)=1,$    coincides now  with a
one sided exit problem (again, since hitting the upper boundary is
impossible). Sine the upper boundary is irrelevant, we will assume
thus w.l.o.g. that $L=\infty$ and $l=0$, or, equivalently, that
  $K<1$ (note that  $\varphi_K(0)=\frac{\lambda+q}{\mu}\left(K-1 \right), \lim_{x\rightarrow \infty}\varphi_K(x)=-\frac{\lambda+q}{\mu}$). 
  
Under the above assumptions, the corresponding Segerdahl's equation becomes integrable (by quadratures) by applying Theorem \ref{FLTU} with $c_0=1$, $c_2=-1$, and $c_1=0$. More specifically, the change of variables given by
\begin{equation}\label{FLchangeDown}
\eta=\sqrt{\frac{-\lambda}{\mu\varphi_{K}(x)}}\bar\eta,
\end{equation}
transforms the Riccati equation in (\ref{FLRicc}) into
\begin{equation}\label{FLQuasiAutDown}
\frac{d \bar\eta}{dx}=\sqrt{\frac{-\lambda\mu}{\varphi_{K}(x)}}(1-\bar\eta^2).
\end{equation}
In order to transform the Riccati equation (\ref{FLQuasiAutDown}) into
an autonomous one, we define the reparametrisation
{\small \begin{equation*}
d\bar x=\sqrt{\frac{-\lambda\mu}{\varphi_{K}(x)}}dx\Longrightarrow
\bar x(x)=\frac
12\sqrt{\frac{\lambda}{q+\lambda}}\log\left(\frac{1-\sqrt{1-K}}{\sqrt{1-K}+1}\frac{\sqrt{1-e^{-2x\mu}K}+1}{1-\sqrt{1-e^{-2x\mu}K}}\right).
\end{equation*}}In consequence, $\bar x(0)=0$ and
$\lim_{x\rightarrow\infty}\bar x(x)=\infty$.
Now, the solution for equation (\ref{FLQuasiAutDown}) reads
$$\bar\eta(\bar x)=\frac{e^{2\bar x}-K_1}{e^{2\bar x}+K_1},$$
with $K_1$ a constant. Hence, in view of the change of variables
(\ref{FLchangeDown}), we obtain that the general solution of the
Riccati equation in (\ref{FLRicc}) under our present assumptions is
$$
\eta(x)=\sqrt{\frac{-\lambda}{\mu\varphi_{K}(x)}}\left(\frac{e^{2\bar x(x)}-K_1}{e^{2\bar x(x)}+K_1}\right)=\sqrt{\frac{\lambda}{q+\lambda}}(1-e^{-2x\mu}K)^{-1/2}\frac{e^{2\bar x(x)}-K_1}{e^{2\bar x(x)}+K_1}.
$$For each value of the function $\eta(x)$, we obtain, in view of the equation $dM/dx=(\mu\eta(x)-\mu)M$,  that the function $M(x)$ reads
\begin{equation*}
M(x)=L\exp\left(\mu\int^x\eta(\bar x)d\bar x-\mu x\right)=\frac
{K_3}2e^{-\mu x}(K_1e^{-\bar x(x)}+e^{ \bar x(x)}).
\end{equation*}
On the other hand, $\Psi(x)=\eta(x)M(x)$. Therefore, we obtain
$$
\Psi(x)=\frac 12
K_3\sqrt{\frac{\lambda}{q+\lambda}}(e^{2x\mu}-K)^{-1/2}(e^{\bar x(x)}-K_1e^{-\bar x(x)}).
$$
Summing up, the solutions for the system (\ref{FLRicc}) in the
particular case admitting the drift (\ref{FLGenSolution}) are
$$\left\{
\begin{aligned}
\Psi(x)&=\frac{K_3}2 \sqrt{\frac{\lambda}{q+\lambda}}(e^{2x\mu}-K)^{-1/2}(e^{\bar x(x)}-K_1e^{-\bar x(x)}),\\
M(x)&=\frac{K_3}2e^{-\mu x}(K_1e^{-\bar x(x)}+e^{ \bar x(x)}).
\end{aligned}\right.
$$
Assume now that $\Psi(0)=M(0)=1$. Therefore,
\begin{equation}\label{FinSol}\left\{
\begin{aligned}
\Psi(x)&=\frac{1}{(1+K_1(K))}\sqrt{\frac{\lambda}{q+\lambda}}(e^{2x\mu}-K)^{-1/2}(e^{\bar x(x)}-K_1(K)e^{-\bar x(x)}),\\
M(x)&=\frac{1}{(1+K_1(K))}e^{-\mu x}(e^{ \bar x(x)}+K_1(K)e^{-\bar x(x)}),
\end{aligned}\right.
\end{equation}
where
$$
K_1(K)=\frac{\sqrt{\frac{\lambda}{q+\lambda}}-\sqrt{1-K}}{\sqrt{\frac{\lambda}{q+\lambda}}+\sqrt{1-K}}.
$$
In order to analyse the behaviour of the above solutions at infinity, it
is necessary to obtain an approximate expression of $\bar x(x)$ at $x>>0$. As $\sqrt{1+y}\simeq 1+y/2$ when
$y\simeq 0$, then
\begin{multline*}
\bar x(x)\simeq\sqrt{\frac{\lambda}{4(q+\lambda)}}\log\left(\frac{\sqrt{1-e^{-2x\mu}K}+1}{|1-\sqrt{1-e^{-2x\mu}K}|}\right)\simeq\\\simeq
\sqrt{\frac{\lambda}{4(q+\lambda)}}\log\left(\frac{4}{e^{-2x\mu}|K|}\right)\simeq
\sqrt{\frac{\lambda}{\lambda+q}}\mu x.
\end{multline*}
Using the above result, it follows that
$$\left\{
\begin{aligned}
\Psi(x)&\simeq\frac{1}{(1+K_1(K))}\sqrt{\frac{\lambda}{q+\lambda}}e^{\mu\left(\sqrt{\frac{\lambda}{\lambda+q}}-1\right)x},\\
M(x)&\simeq\frac{1}{(1+K_1(K))}e^{\mu\left(\sqrt{\frac{\lambda}{\lambda+q}}-1\right) x}.
\end{aligned}\right.
$$
Consequently, it turns out that $\lim_{x\rightarrow
\infty}\Psi(x)=\lim_{x\rightarrow \infty}M(x)=0$.

In order to illustrate the behaviour of the ruin probabilities (\ref{FinSol}), we provide the following figure showing a particular case of our solutions and its corresponding drift $\varphi_K(x)$.

\begin{figure}[hbt]
\begin{center}
\includegraphics[scale=0.7]{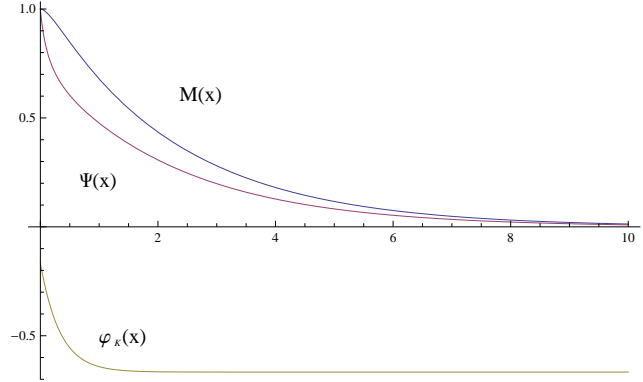}
\end{center}
\caption{Ruin probabilities and drift for $\mu=1.5, \lambda=q=1/2$, and $K=0.75$.}
\label{f1}
\end{figure}

\section{Acknowldgements}
This work was partially supported by research projects MTM2009-11154, E24/1 (DGA), and the CNRS-LEA Math-Mode funding 2010.

\end{document}